\documentclass[12pt]{article}
\usepackage{amsmath,amsthm,amssymb}
\usepackage{mathtools}

\newtheorem{theorem}{Theorem}[section] 
\newtheorem{lemma}{Lemma}[section]

\newtheorem{remark}{Remark}[section]

\newtheorem{ass}{Assumption}[section]

\newtheorem{proposition}{Proposition}[section]

\numberwithin{equation}{section}

\newcommand{\D}{{\rm d}}
\newcommand{\dx}{\, \D x}

\newcommand{\dis}{\displaystyle}

\newcommand{\rz}{\mathbb{R}}
\newcommand{\nz}{\mathbb{N}}

\newcommand{\eps}{\varepsilon}
\newcommand{\iom}{\int_{\Omega}}
\newcommand{\ib}{\int_{B}}

\newcommand{\arsp}{\hspace*{\arraycolsep}}

\title{Splitting-type variational problems with asymmetrical growth conditions}

\author{Michael Bildhauer \& Martin Fuchs}
\date{}

\newcommand{\reff}[1]{(\ref{#1})}

\usepackage{hyperref}

\hypersetup{
   pdfnewwindow=true,
   colorlinks=false,
   linkcolor=black,
   citecolor=green,
   filecolor=magenta,
   urlcolor=cyan,
}

\begin{document}

\parindent0em
\maketitle

\newcommand{\op}[1]{\operatorname{#1}}
\newcommand{\bv}{\op{BV}}
\newcommand{\mub}{\overline{\mu}}
\newcommand{\muhat}{\hat{\mu}}

\newcommand{\hypref}[2]{\hyperref[#2]{#1 \ref*{#2}}}
\newcommand{\hypreff}[1]{\hyperref[#1]{(\ref*{#1})}}

\newcommand{\ob}[1]{^{(#1)}}

\newcommand{\xh}{\Xi}
\newcommand{\oh}[1]{O\left(#1\right)}
\newcommand{\xn}{w_{1}}
\newcommand{\yn}{w_{2}}
\newcommand{\On}{\hat{\Omega}}
\newcommand{\uq}{\hat{u}}
\newcommand{\gameps}{\Gamma_{Q}}
\newcommand{\se}{{s_Q}}
\begin{abstract}{}
Splitting-type variational problems
\[
\iom \sum_{i=1}^n f_i(\partial_i w) \dx \to \min 
\]
with superlinear growth conditions are studied by assuming 
\[
h_i(t) \leq  f''_i(t)  \leq H_i(t) 
\]
with suitable functions $h_i$, $H_i$: $\rz \to \rz^+$, $i=1$, \dots , $n$, 
measuring the growth and ellipticity of the energy density. Here, as the main feature,
a symmetric behaviour like $h_i(t)\approx h_i(-t)$ and $H_i(t) \approx H_i(-t)$ for large $|t|$
is not supposed.

Assuming quite weak hypotheses as above, we establish higher integrability 
of $|\nabla u|$ for local minimizers $u\in L^\infty(\Omega)$ by
using a Caccioppoli-type inequality with some power weights of negative
exponent.
\end{abstract}

\newcommand{\platz}{\vspace{2ex}}
\newcommand{\absatz}{
\platz
\centerline{\rule[1ex]{3cm}{0.01cm}}
\platz
}

\newcommand{\os}{\overline{s}_\gamma}
\newcommand{\oslimit}{\overline{s}}
\newcommand{\ls}{\varkappa}

\newcommand{\gamf}{{\theta_{i}}}
\newcommand{\gamg}{{\gamma_{i}}}

\parindent0ex
\section{Introduction}\label{in}


Suppose that $\Omega\subset \rz^n$ is a bounded Lipschitz domain and 
consider the variational integral
\begin{equation}\label{in 1}
J[w] := \iom f\big(\nabla w\big) \dx 
\end{equation}
of splitting-type, i.e.
\begin{equation}\label{in 2}
f:\, \rz^n \to \rz\, , \quad f(Z) = \sum_{i=1}^n f_i(Z_i)
\end{equation}
with strictly convex functions $f_i$: $\rz \to \rz$ of class $C^2(\rz)$, $i=1$, \dots , $n$, satisfying in addition some
suitable superlinear growth and ellipticity conditions. \\

Problem \reff{in 1}, \reff{in 2} serves as a prototype for non-uniformly elliptic variational problems and it is well
known that the ratio of the highest and the lowest eigenvalue of $D^2f$ is the crucial quantity for proving the regularity 
of solutions (see, e.g., \cite{Gi:1987_1}, \cite{Ma:1989_1}, \cite{Ma:1996_1}). The reader will find an extensive overview including 
different settings of non-uniformly elliptic variational problems in the recent paper \cite{BM:2020_1}.
Without going into further details we refer to the series of references given in this paper.\\

In Section 1.3 of \cite{BM:2020_1}, the authors consider general growth conditions which, roughly speaking, means that
the energy density $f$ is controlled in the sense of
\begin{equation}\label{in 3}
g(|Z|)|\xi|^2 \leq  D^2 f(Z)(\xi,\xi)\, , \quad |D^2f(Z)| \leq G(|Z|)\, ,
\end{equation}
with suitable functions $g(t)$, $G(t)$: $\rz^+_0 \to \rz^+$. Then, under appropriate assumptions on $g$, $G$, a
general approach to regularity theory is given in \cite{BM:2020_1}.\\

Our note is motivated by the observation, that in \reff{in 2} there is no obvious reason to assume some kind of 
symmetry for the functions $f_i$, i.e.~in general we have $f_i(t) \not= f_i(-t)$ and, as one model case, we just 
consider ($q_i^{\pm} \geq 1$, $i=1,\dots , n$)
\begin{equation}\label{in 4}
f_i(t) \approx |t|^{q_i^-}\quad\mbox{if}\quad t \ll -1\, , \qquad
f_i(t) \approx |t|^{q_i^+}\quad\mbox{if}\quad t \gg 1 \, .
\end{equation}
Then, both for $t\ll 1$ and for $t \gg 1$, the functions $f_i$ just behave like a uniform  power of $|t|$. Nevertheless, 
the power $q_i^-$ enters the left-hand side of \reff{in 3} and $q^+_i$ is needed on the right-hand side of \reff{in 3}.\\ 

This motivates to study the model case \reff{in 2} and to establish regularity results for solutions under the weaker assumption
\begin{equation}\label{in 5}
h_i(t) \leq f_i''(t) \leq H_i(t) \quad t \in \rz\, ,
\end{equation}
with suitable functions $h_i$, $H_i$: $\rz \to \rz^+$, $i=1,\dots n$.\\

There is another quite subtle difficulty in studying regularity of solutions to splitting-type variational problems: in
\cite{BFZ:2007_1} the authors consider variational integrals of the form ($1 \leq  k < n$)
\begin{equation}\label{in 6}
I[w,\Omega]  = \iom \Big[f(\partial_1 w,\dots ,\partial_{k} w) + g(\partial_{k+1}w,\dots , \partial_n w)\Big] \dx \, , 
\end{equation}
where $f$ and $g$ are of $p$ and $q$-growth, respectively ($p$, $q >1$). Then the regularity of bounded solutions follows
in the sense of \cite{BFZ:2007_1}, Theorem 1.1, without any further condition relating $p$ and $q$. The proof argues
step by step and works since the energy density splits into two parts. If, as supposed in \reff{in 2}, the energy density
splits in more than two components, then one has to be more careful dealing with the exponents and some
more restrictive (but still quite weak) assumptions have to be made. In this sense Remark 1.3 of \cite{BFZ:2007_1}
might be a little bit misleading. We note that a splitting structure into two components as supposed in \reff{in 6}
is also assumed, e.g., in \cite{Br:2010_1} and related papers.\\

In the following we consider the variational integral \reff{in 1}, \reff{in 2} defined on the energy class
\[
E_f(\Omega) := \Bigg\{w\in W^{1,1}(\Omega):\, \iom f(\nabla w) \dx < \infty \Bigg\}\, .
\]
We are interested in local minimizers $u$: $\Omega \to \rz$ of class $E_f(\Omega)$, i.e.~it
holds that
\begin{equation}\label{in 7}
\iom f(\nabla u) \dx \leq \iom f(\nabla w) \dx
\end{equation}
for all $w\in E_f(\Omega)$ such that $\op{spt}(u-w) \Subset \Omega$.\\

{\bf Notation.} We will always denote by $q_i^+ >1$, $q_i^- >1$, $1\leq i \leq n$, real exponents and we let
for fixed $1 \leq i \leq n$
\begin{equation}\label{in 8}
\underline{q}_i := \min\{q_i^{\pm}\}\, ,\quad\overline{q}_i := \max\{q_i^{\pm}\} \, .
\end{equation}
Moreover, we let
\[
\Gamma:\; [0,\infty) \to \rz\, ,\quad \Gamma(t) = 1 + t^2\, .
\]
Recalling the idea sketched in \reff{in 4}, \reff{in 5} we denote by
$h_i$ and $H_i$, $i=1$, \dots , $n$, functions $\rz \to \rz^+$ such that with positive constants
$\underline{a}_i$, $\overline{a}_i$
\begin{equation}\label{in 9}
\left. \begin{array}{l} 
\dis\underline{a}_i \Gamma^{\frac{q_i^{-}-2}{2}}(|t|)\quad\mbox{if $t < -1$}\\[2ex]
\dis \underline{a}_i \Gamma^{\frac{q_i^{+}-2}{2}}(|t|)\quad\mbox{if $t > 1$}\\
\end{array}\right\}
\leq  h_i(t)
\end{equation}
and 
\begin{equation}\label{in 10}
 H_i(t) 
\leq 
\left\{ \begin{array}{l} 
\dis\overline{a}_i \Gamma^{\frac{q_i^{-}-2}{2}}(|t|)\quad \mbox{if $t < -1$}\\[2ex]
\dis \overline{a}_i \Gamma^{\frac{q_i^{+}-2}{2}}(|t|)\quad \mbox{if $t > 1$}\\
\end{array}\right.  \, .
\end{equation}

\vspace*{2ex}
As a general assumption we consider functions $f_i$: $\rz \to [0,\infty)$ of class $C^2(\rz)$, $i=1$, \dots , $n$,
such that  for all $t\in \rz$
\begin{equation}\label{in 11}
h_i(t) \leq  f''_i(t)  \leq H_i(t) 
\end{equation}
and note that \reff{in 11} immediately implies for all $i \in \{1, \dots , n\}$ with constants $b_i>0$
\begin{equation}\label{in 12}
|f_i'(t)|  \leq  b_i 
\left\{ \begin{array}{l} 
\dis \Gamma^{\frac{q_i^{-}-1}{2}}(|t|)\;\mbox{if $t < -1$}\\[1ex]
\dis \Gamma^{\frac{q_i^{+}-1}{2}}(|t|)\;\mbox{if $t > 1$}\\
\end{array}\right\}\, .
\end{equation}
Moreover we obtain for all $i=1$, \dots , $n$ with constants $\underline{c}_i$, $\overline{c}_i > 0$
\begin{equation}\label{in 13}
\underline{c}_i \left\{ \begin{array}{l} 
\dis \Gamma^{\frac{q_i^{-}}{2}}(|t|)\;\mbox{if $t < -1$}\\[2ex]
\dis \Gamma^{\frac{q_i^{+}}{2}}(|t|)\;\mbox{if $t > 1$}\\
\end{array}\right\}
\leq  f_i(t) 
\leq \overline{c}_i
 \left\{ \begin{array}{rcl} 
\dis \Gamma^{\frac{q_i^{-}}{2}}(|t|);\mbox{if $t < -1$}\\[2ex]
\dis \Gamma^{\frac{q_i^{+}}{2}}(|t|)\;\mbox{if $t > 1$}\\
\end{array}\right\} \, .
\end{equation}

With this notation our main result reads as follows.

\begin{theorem}\label{main}
Suppose that for $i=1$, \dots , $n$ the functions $f_i$: $\rz \to [0,\infty)$ are of class $C^2(\rz)$ and satisfy \reff {in 11}
with $h_i$, $H_i$ given in \reff{in 9}, \reff{in 10}. \\

With the notation \reff{in 8} we assume in addition that we have for every fixed $1 \leq i \leq n$
\begin{eqnarray}
\label{in 14}
\overline{q}_j &<& 2 \underline{q}_i + 2   \qquad \mbox{for all}\qquad i <  j \leq n\, ,\\[1ex]
\label{in 15}
\overline{q}_j &< &3 \underline{q}_i +2 \qquad\mbox{for all}\qquad 1 \leq j < i \, .
\end{eqnarray}
If $u \in L^\infty(\Omega) \cap E_f(\Omega)$ denotes a local minimizer of \reff{in 1}, \reff{in 2}, i.e.~of
\[
J[w] = \iom \Bigg[\sum_{i=1}^n f_i(\partial_i w)\Bigg]\dx \, , 
\]
then there exists a real number $\delta > -1/2$ such that for every $1\leq i \leq n$
\begin{equation}\label{in 16}
\ib f_i (\partial_i u) \Gamma^{1+\delta}(|\partial_i u|)  \eta^{2k}\dx \leq c \, .
\end{equation}
\end{theorem}

\begin{remark}\label{in rem 2}
As outlined in Remark \ref{proof rem 2} and Remark \ref{proof rem 3} below, we recover the results of  \cite{BFZ:2007_1}
in the sense that \reff{in 14} and \reff{in 15} are superfluous in the case $n=2$ (or related situations) and
\[
q_1^+=q_1^-\, , \quad q_2^+=q_2^- \, .
\]
\end{remark}

\vspace*{2ex}

Theorem \ref{main} describes the typical situation we have in mind. The proof however is not limited to this particular case
which leads to the generalized version stated in Theorem \ref{main 2} below. \\

In Section \ref{regular} we shortly sketch a regularization procedure via Hilbert-Haar solutions while Section \ref{proof}
presents the main inequalities for the iteration procedure of Section \ref{iter}. This completes the proof ot Theorem \ref{main 2}
and hence Theorem \ref{main}.

\section{Precise assumptions on $f$}\label{pre}
The suitable larger class of admissible energy densities is given by the following assumption.

\begin{ass}\label{pre ass 1}
The energy density $f$,
\[
f:\, \rz^n \to \rz\, , \quad f(Z) = \sum_{i=1}^n f_i(Z_i) \, ,
\]
introduced in \reff{in 2} is supposed to satisfy the following hypotheses.
\begin{enumerate}
\item The function $f_i$: $\rz \to [0,\infty)$, $i=1$, \dots , $n$, is of class $C^2(\rz)$ and 
for all $t \in \rz$ we have $f_i''(t) >0$.

 For $1 \leq i \leq n$ we suppose superlinear growth in the sense of 
\[ 
 \lim_{t\to \pm \infty} |f_i'(t)| = \infty
\]
and at most of polynomial growth in the sense that for some $s> 0$ we have for $|t|$ sufficiently large
\[
f_(t) \leq c |t|^s\quad \mbox{with a finte constant $c$}\, .
\] 

\item For $i\in \{1,\dots ,n\}$ with exponents $\delta_i \geq 0$, $\theta_i  \geq 0$ satisfying
\begin{equation}\label{pq}
\theta_i < 1-\delta_i
\end{equation}
we suppose that for all $|t|$ sufficiently large
\begin{eqnarray}
\label{pre 1}
c_1 \Gamma^{1-\delta_i}(|t|) f''_i (t) &\leq& f_i (t) \arsp \leq\arsp
c_2 f''_i(t) \Gamma^{1+ \theta_i}(|t|)\, ,\\[2ex]
\label{pre 2}
|f'_i(t)|^2 &\leq& c_3 f''_i(t) f_i(t) \Gamma^{\theta_i}(|t|)\, ,
\end{eqnarray}
where $c_1$, $c_2$ and $c_3$ denote positive constants.
\item We let
\[
\Gamma^{\frac{q_i^{\pm}}{2}}(t) =  
\left\{ \begin{array}{l} 
\dis \Gamma^{\frac{q_i^{-}}{2}}(|t|)\;\mbox{if $t < 0$}\\[2ex]
\dis \Gamma^{\frac{q_i^{+}}{2}}(|t|)\;\mbox{if $t \geq 0$}
\end{array}\right\} \, .
\]
and suppose that $f_i$, $i=1$, \dots ,$n$, satisfies with positive constants $c_4$, $c_5$ and for $|t|$
sufficiently large
\begin{equation}\label{pre 3}
c_4 \Gamma^{\frac{q_i^\pm}{2}}(|t|) \leq  f_i(t) \leq c_5  \Gamma^{\frac{q_i^\pm}{2}}(|t|) \, .
\end{equation}
\end{enumerate}
\end{ass}

\begin{remark}\label{pre rem 1}
\begin{enumerate}
\item If $f_i$ is a power growth function like, e.g., $f_i(t) = (1+t^2)^{p_i/2}$, $p_i > 1$ fixed, then we have
\[
c \Gamma(|t|) f_i''(t) \leq f_i(t) \leq c\Gamma(|t|) f_i''(t) \, ,
\]
i.e.~\reff{pre 1}.  The same is true for our asymmetric model case given by \reff{in 9} -- \reff{in 13}.

\item By convexity it is well known (see, e.g., \cite{Da:2015_1}, exercise 1.5.9, p.~53) that the right-hand side of \reff{pre 2} 
and the right-hand side of \reff{in 13} imply \reff{pre 3}.

\item The condition \reff{pq} formally corresponds with the condition $q < p+2$ in the standard $(p,q)$-case
(see, e.g., \cite{Bi:1818}, Chapter 5, and the references quoted therein).

\item Assumption \ref{pre ass 1}, $iii$) is assumed w.l.o.g. 
In fact, on account of $f_i''>0$ we know that for any $i\in \{1,\dots , n\}$ the function $f'_i$ is an increasing function, 
and by Assumption \ref{pre ass 1}, $i$), we let
\[
s^+ := \inf_{s} \lim_{t\to \infty} \frac{f'(t)}{t^s} < \infty \, , 
\quad s^- := \inf_{s} \lim_{t\to -\infty} \frac{|f'(t)|}{|t|^s} < \infty \, .
\] 
Then for arbitrary small $\eps >0$ we have the right-hand side of \reff{pre 3} with exponent $s^\pm+1+\eps$
and the left-hand side with exponent $s^\pm + 1 - \eps$. Going through the proof of Theorem \ref{main 2} we may
suppose \reff{pre 3} w.l.o.g.
\end{enumerate}
\end{remark}

\begin{theorem}\label{main 2}
Suppose that we have Assumption \ref{pre ass 1}.  With the above notation we assume in addition that
we have for every fixed $1 \leq i \leq n$
\begin{eqnarray}
\label{pre 4}
\overline{q}_j &<& \frac{2 \underline{q}_i (1-\delta_i)}{1+2 \theta_i} + 2 (1-\delta_i)  
\qquad \mbox{for all}\qquad i <  j \leq n\, ,\\[1ex]
\label{pre 5}
\overline{q}_j &< &\frac{2}{1+2\theta_i}\Bigg[\frac{\underline{q}_i}{2}(1-\delta_i)\Big[2+\frac{1}{1- \delta_i}\Big] 
- \theta_i (1+\overline{q}_j)\Bigg]+2 (1-\delta_i)\nonumber \\[1ex]
&& \qquad\mbox{for all}\qquad 1 \leq j < i \, .
\end{eqnarray}
If $u \in L^\infty(\Omega) \cap E_f(\Omega)$ denotes a local minimizer of \reff{in 1}, \reff{in 2}, i.e.~of
\[
J[w] = \iom \Bigg[\sum_{i=1}^n f_i(\partial_i w)\Bigg]\dx \, , 
\]
then there exists a real number $\delta > -1/2$ such that for every $1\leq i \leq n$
\begin{equation}\label{pre 6}
\ib f_i (\partial_i u) \Gamma^{1+\delta}(|\partial_i u|)  \eta^{2k}\dx \leq c \, .
\end{equation}
\end{theorem}

\begin{remark}\label{pre rem 2}
In particular we note that \reff{pre 4}, \reff{pre 5} reduce to \reff{in 14}, \reff{in 15} 
for $\delta_i$, $\theta_i$ sufficiently small.
\end{remark}

\newcommand{\vart}{\vartheta_{i}}

\section{Some remarks on regularization}\label{regular}

We have to start with a regularization procedure such that the expressions given below are well defined.
We follow Section 2 of \cite{BFZ:2007_1} and fix a ball $D  \Subset \Omega$. If $u$ denotes the local
minimizer under in the sense of \reff{in 7} and if $\eps > 0$ is sufficiently small, we consider the mollification $(u)_\eps$
of $u$ w.r.t.~the radius $\eps$.  We consider the Dirichlet-problem
\[
\int_D \sum_{i=1}^n f_i(\partial_i w) \dx \to \min 
\]
among all Lipschitz mappings $\overline{D}\to \rz$ with boundary data $(u)_\eps$. According to, e.g., \cite{MM:1984_1},
there exits a unique (Hilbert-Haar) solution $u_\eps$ to this problem.\\

Exactly as outlined in \cite{BFZ:2007_1} Lemma 2.1 and Lemma 2.2 we obtain:

\begin{lemma}\label{regular lem 1}
Let $\underline{q} := \min_{1\leq i \leq n} \underline{q}_i$
\begin{enumerate}
\item We have as $\eps \to 0$
\[
u_\eps \rightharpoondown u \quad\mbox{in}\, W^{1,\underline{q}}(D)\, , \qquad
\int_D \sum_{i=1}^n f_i(\partial_i u_\eps) \dx \to \int_D \sum_{i=1}^n f_i(\partial_i u) \dx \, .
\]
\item There is a finite constant $c >0$ not depending on $\eps$ such that 
\[
\| u_\eps\|_{L^{\infty}(D)} \leq c \, .
\]
\item For any $\alpha < 1$ we have $u_\eps \in C^{1,\alpha}(D)\cap W^{2,2}_{\op{loc}}(D)$.
\end{enumerate}
\end{lemma}

We then argue as follows: consider a local minimizer $u$ of \reff{in 1}, \reff{in 2} and the approximating sequence
$\{u_\eps\}$ minimizing the functional
\begin{equation}\label{regular 1}
J[w,D] := \int_D \sum_{i=1}^n f_i(\partial_i w_i) \dx
\end{equation}
w.r.t.~the data $(u)_\eps$. In particular we have a sequence of local $J[w,D]$-minimizers. We apply the a priori
results of the next section to $u_\eps$ and Theorem \ref{main} follows from Lemma \ref{regular lem 1}
passing to the limit $\eps \to 0$.

\section{General inequalities}\label{proof}

The main result of this section is Proposition \ref{proof prop 2} which is not depending on the particular
structure \reff{in 9}.\\

We will rely on the following variant of Caccioppoli's inequality which was first
introduced in \cite{BF:2020_3}. We also refer to Section 6 of \cite{BF:2022_2} on Caccioppoli-type inequalities
involving powers with negative exponents, in particular we refer to Proposition 6.1.

\renewcommand{\iom}{\int_{D}}

\begin{lemma}\label{proof lem 1}
Fix $l\in \nz$ and suppose that $\eta \in C^\infty_0(D)$, $0 \leq \eta \leq 1$.  If we consider a local minimizer
$u \in W^{1,\infty}_{\op{loc}}(D) \cap W^{2,2}_{\op{loc}}(D)$ of the variational functional
\[
I[w] = \iom g(\nabla w) \dx
\]
with energy density $g$: $\rz^n \to \rz$ of class $C^2$ satisfying $D^2 g(Z)(Y,Y) > 0$ for all $Y$, $Z\in \rz^n$, then
for any fixed $i \in \{1,\dots , n\}$ we have
\begin{eqnarray*}
\lefteqn{\iom D^2 g (\nabla u)\big(\nabla \partial_i u, \nabla \partial_i u\big) \eta^{2l} 
\Gamma^{\beta}(|\partial_i u|) \dx}\nonumber \\ 
&& \leq c \iom D^2 g (\nabla u) (\nabla \eta,\nabla \eta)\eta^{2l-2} \Gamma^{1+\beta}(|\partial_i u|) \dx
\end{eqnarray*}
for any $\beta > - 1/2$.
\end{lemma}

To the end of our note we always consider a fixed ball $B=B_{2r}(x_0) \Subset D$. With this notation
we have the following auxiliary proposition. 

\begin{proposition}\label{proof prop 1}
Suppose that we have $i)$ of Assumption \ref{pre ass 1} and let $\eta \in C^{\infty}_0(B)$, $0 \leq \eta \leq 1$, 
$\eta \equiv 1$ on $B_r(x_0)$, $|\nabla \eta| \leq c/r$. 
Moreover, we assume that  $u \in L^\infty(D) \cap W^{1,\infty}_{\op{loc}}(D) \cap W^{2,2}_{\op{loc}}(D)$.\\

Then we have for fixed $\gamma \in \rz$, for all $k>0$ sufficiently large and for $i=1$, \dots, $n$ the starting inequalities (no summation w.r.t.~$i$)
\begin{eqnarray}\label{proof 1}
\lefteqn{\ib f_i(\partial_i u) \Gamma^{1+\gamma}(|\partial_i u|) \eta^{2k}\dx}\nonumber\\[2ex]
&\leq & c \Bigg[1 + \ib |\partial_i \partial_i u| \Gamma^{\gamma}(|\partial_i u|) f_i(\partial_i u) \eta^{2k}\dx\nonumber\\[2ex]
&&+ \ib |\partial_i\partial_i u| \, |f_i'|(\partial_i u) \, \Gamma^{\frac{1}{2}+\gamma}(|\partial_i u|) \eta^{2k} \dx\Bigg] 
\end{eqnarray}
\end{proposition}

\begin{remark}\label{proof rem 1} 
\begin{enumerate}
\item The idea of the proof of Proposition \ref{proof prop 1} is based on an 
integration by parts using the boundedness of $u$. 
An Ansatz of this kind was already made by  Choe \cite{Ch:1992_1}, where all relevant quantities are depending on $|\nabla u|$.
Here the main new feature is to work with the energy density $f$ which is not depending on the modulus of $\nabla u$. 
\item We note that for the proof of Proposition \ref{proof prop 1} no minimizing property of $u$ is needed.
\end{enumerate}
\end{remark}

{\it Proof of Proposition \ref{proof prop 1}.} With $i\in \{1,\dots , n\}$ fixed we obtain using an integration by parts
\begin{eqnarray}\label{proof 2}
\lefteqn{\ib f_i(\partial_i u) \Gamma^{1+\gamma}(|\partial_i u|) \eta^{2k}\dx}\nonumber\\[2ex] 
&=& \ib |\partial_i u|^2 f_i(\partial_i u) \Gamma^{\gamma}(|\partial_i u|)\eta^{2k}\dx + 
\ib f_i(\partial_i u) \Gamma^{\gamma}(|\partial_i u|)\eta^{2k}\dx\nonumber\\[2ex]
&=& 
- \ib u \partial_i \Big[\partial_i u f_i(\partial_i u) \Gamma^{\gamma}(|\partial_i u|) \eta^{2k}\Big]\dx
+ \ib f_i(\partial_i u) \Gamma^{\gamma}(|\partial_i u|)\eta^{2k}\dx
\nonumber\\[2ex]
&\leq & c \ib |\partial_i \partial_i u| \Gamma^{\gamma}(|\partial_i u|) f_i(\partial_i u) \eta^{2k} \dx\nonumber\\[2ex]
&&+ c \ib |\partial_i \partial_i u|\, |\partial_i u| \, |f'_i|(\partial_i u)\, \Gamma^{\gamma}(|\partial_i u|) \eta^{2k} \dx
\nonumber\\[2ex]
&&+ c \ib |\partial_i u| f_i(\partial_i u) \Gamma^{\gamma}(|\partial_i u|) \eta^{2k-1} |\partial_i \eta| \dx 
+ \ib f_i(\partial_i u) \Gamma^{\gamma}(|\partial_i u|) \eta^{2k}\dx
\nonumber\\[2ex]
&=& I_{1,i} + I_{2,i} + I_{3,i} + I_{4,i}\, .
\end{eqnarray}
In \reff{proof 2} we discuss $I_{3,i}$: for $\eps >0$ sufficiently small we estimate
\begin{eqnarray}\label{proof 3}
I_{3,i} & \leq&  \ib |\partial_i u| f_i^{\frac{1}{2}}(\partial_i u) \Gamma^{\frac{\gamma}{2}}(|\partial_i u|) \eta^{k} 
f_i^{\frac{1}{2}}(\partial_i u) \Gamma^{\frac{\gamma}{2}}(|\partial_i u|)\eta^{k-1} |\nabla \eta| \dx\nonumber\\[2ex]
&\leq & \eps \ib |\partial_i u|^2 f_i(\partial_i u) \Gamma^{\gamma}(|\partial_i u|) \eta^{2k} dx\nonumber\\[2ex]
&&+c(\eps,r) \ib f_i(\partial_i u) \Gamma^{\gamma}(|\partial_i u|) \eta^{2k-2}\dx \, .
\end{eqnarray}
The first integral on the right-hand side of \reff{proof 3} is absorbed in the left-hand side of \reff{proof 2}, i.e.
\begin{eqnarray}\label{proof 4}
\lefteqn{\ib f_i (\partial_i u) \Gamma^{1+\gamma}(|\partial_i u|) \eta^{2k}\dx}\nonumber\\[2ex]
&\leq &  I_{1,i} +  I_{2,i}
+ c(\eps,r)  \ib f_i(\partial_i u)  \Gamma^{\gamma}(|\partial_i u|)\eta^{2k-2}\dx\nonumber\\[2ex] 
&&+ \ib f_i (\partial_i u) \Gamma^{\gamma}(|\partial_i u|)\eta^{2k}\dx\nonumber\\[2ex]
&\leq &  I_{1,i} +  I_{2,i} + c(\eps,r)  \ib f_i(\partial_i u)  \Gamma^{\gamma}(|\partial_i u|)\eta^{2k-2}\dx \, .
\end{eqnarray}

Discussing the remaining integral we recall that the function $f_i(t) \Gamma^{1+\gamma}(|t|)$ is at most of polynomial growth, 
hence we may apply the auxiliary Lemma \ref{proof lem 2} below to the functions 
$\varphi(t) = f_i(t) \Gamma^{\gamma}(|t|)$ and 
$\psi(t) := f_i(t) \Gamma^{1+\gamma}(|t|)$ with the result that for some $\rho >0$ and for all $t \in \rz$
\begin{equation}\label{proof 5}
f_i(t) \Gamma^{\gamma}(|t|) \leq 
c \big[f_i(t) \Gamma^{1+\gamma}(|t|)\big]^{\frac{1}{\rho}} + c \, .
\end{equation}

With \reff{proof 5} we estimate for $\tilde{\eps} >0$ sufficiently small and for $k > \rho^* = \rho/(\rho-1)$
\begin{eqnarray}\label{proof 6}
\lefteqn{c(\eps,r) \ib f_i(\partial_i u)  \Gamma^{\gamma}(|\partial_i u|) \eta^{2k-2}\dx}\nonumber\\[2ex]
& \leq &
c(\eps,r) \ib  \big[f_i(\partial_i u)  \Gamma^{1+\gamma}(|\partial_i u|)\big]^{\frac{1}{\rho}}
\eta^{\frac{2k}{\rho}} \eta^{\frac{2k}{\rho^*}-2}\dx +c\nonumber\\[2ex]
& \leq & \tilde{\eps} \ib f_i (\partial_i u) \Gamma^{1+\gamma}(|\partial_i u|)\eta^{2k}\dx
+ c(\tilde{\eps},\eps,r) \ib \eta^{2(k-\rho^*)} \dx +c \, .
\end{eqnarray}
The inequalities \reff{proof 4} and \reff{proof 6} complete the proof of the proposition by absorbing
the first integral on the right-hand side of \reff{proof 6} in the left-hand side of \reff{proof 4}.  \qed \\

It remains to give an elementary proof of the following auxiliary Lemma.
\begin{lemma}\label{proof lem 2}
For $m\in \nz$ we consider functions $\varphi$, $\psi$: $\rz^m \to [0,\infty)$ such that $\psi(X) \leq c \Gamma^\tau(|X|)$ 
for some $\tau >0$ and for all $X\in \rz^m$.
Suppose that we have for some $\eps >0$ and for all $X \in \rz^n$
\[
\varphi(X) \leq c \Gamma^{-\eps}(|X|) \psi(X)\, .
\]
Then there exists a real number $\rho > 1$ and a constant $C >0$ such that
\[
\varphi(X) \leq \big[\psi(X)\big]^{\frac{1}{\rho}} + C\, .
\]
\end{lemma} 
{\it Proof.} Let $\delta := \eps/\tau$, i.e.~for all $X \in \rz^m$
\[
1+ \psi^\delta  \leq 1+ \Gamma^\eps \leq 2 \Gamma^\eps\, ,
\]
hence we have by assumption
\begin{eqnarray*}
\varphi(X) &\leq & c \big[ 1 + \psi^{\delta}(X)\big]^{-1} \psi(X)\\[2ex]
& \leq & \left\{\begin{array}{ccl}
c& \mbox{if} & \psi^\delta(X) \leq 1\\[2ex]
c \psi^{1-\delta}(X) &\mbox{if}& \psi^{\delta}(X) > 1
\end{array}\right\}\, .
\end{eqnarray*}
The lemma follows with the choice $\rho = 1/(1- \delta)$. \qed\\

With the help of Proposition \ref{proof prop 1} we now establish the main inequality of this section.

\begin{proposition}\label{proof prop 2}
Suppose that we have Assumption \ref{pre ass 1} and let $\eta \in C^{\infty}_0(B)$, $0 \leq \eta \leq 1$, 
$\eta \equiv 1$ on $B_r(x_0)$, $|\nabla \eta| \leq c/r$. 
Moreover, we assume that  $u \in L^\infty(D) \cap W^{1,\infty}_{\op{loc}}(D) \cap W^{2,2}_{\op{loc}}(D)$  
is a local minimizer of \reff{regular 1}.\\

For $i\in \{1,\dots ,n\}$ we choose $\eps_i$ satisfying $\theta_i < \eps_i < 1- \delta_i$ (recall \reff{pq})
and let $\gamma_i+ \eps_i =: \beta_i$, where we always suppose in the following that $\beta_i > -1/2$.\\

Then we have for any sufficiently large real number $k>0$
\begin{equation}\label{proof 7}
\ib f_i (\partial_i u) \Gamma^{1+\gamma_{i}}(|\partial_i u|)  \eta^{2k}\dx \leq c
\sum_{j\not= i}  \ib f''_{j}(\partial_j u) \Gamma^{1+\beta_i}(|\partial_i u|) \eta^{2k-2}\dx \, .
\end{equation}
\end{proposition}

{\it Proof.} We recall the starting inequality \reff{proof 1},
\begin{equation}\label{proof 8}
\ib f_i (\nabla u) \Gamma^{1+\gamma_{i}}(|\partial_i u|) \eta^{2k}\dx \leq c \Big[1+ I_{1,i} + I_{2,i}\Big] \, ,
\end{equation}
where we fix $i \in \{1,\dots , n\}$. We estimate for fixed $\beta_i$ as above 
\begin{eqnarray}\label{proof 9}
I_{1,i} &=& \ib |\partial_i \partial_i u| {f''}_i^{\frac{1}{2}}(\partial_i u) \Gamma^{\frac{\beta_i}{2}}(|\partial_i u|)
(f''_i)^{-\frac{1}{2}}(\partial_i u)  \Gamma^{-\frac{\beta_i}{2}}(|\partial_i u|)\nonumber\\[2ex]
&&\qquad \cdot  \Gamma^{\gamma_{i}}(|\partial_i u|)f_i(\partial_i u)\eta^{2k}\dx\nonumber\\[2ex]
&\leq & c \ib f''_i(\partial_i u) |\partial_i\partial_i u|^2 \Gamma^{\beta_i}(|\partial_i u|) \eta^{2k}\dx\nonumber\\[2ex] 
&&+ c \ib (f''_i)^{-1}(\partial_i u) \Gamma^{\gamma_{i}-\eps_i}(|\partial_i u|) f_i^2(\partial_i u) \eta^{2k}\dx \, .
\end{eqnarray}

The second integral on the right-hand side of \reff{proof 9} is handled with the help of the right-hand side of \reff{pre 1} 
using in addition Lemma \ref{proof lem 2} (recalling $\eps_i > \theta_i$)
\begin{eqnarray}\label{proof 10}
\lefteqn{\ib (f''_i)^{-1}(\partial_i u) \Gamma^{\gamma_{i}-\eps_i}(|\partial_i u|) f_i^{2}(\partial_i u) \eta^{2k}\dx}\nonumber\\[2ex]
 & \leq &\ib \Big[f_i(\partial_i u) \Gamma^{1+\gamma_{i} - (\eps_i-\theta_i)}(|\partial_i u|) \Big] \eta^{2k}\dx\nonumber\\[2ex]
&\leq &
\ib \Big[f_i(\partial_i u) \Gamma^{1+\gamma_{i}}(|\partial_i u|) \Big]^{\frac{1}{\rho}}\eta^{\frac{2k}{\rho}}
\eta^{\frac{2k}{\rho*}}\dx + c \nonumber\\[2ex]
 &\leq & \eps \ib f_i(\partial_i u) \Gamma^{1+\gamma_{i}}(|\partial_i u|)\eta^{2k}\dx + c(\eps,r) \, .
\end{eqnarray}

Absorbing terms it is shown up to now (using \reff{proof 8} - \reff{proof 10})
\begin{eqnarray}\label{proof 11}
\lefteqn{\ib f_i (\partial_i u) \Gamma^{1+\gamma_{i}}(|\partial_i u|) \eta^{2k}\dx}\nonumber\\[2ex] 
&\leq & 
c \Bigg[1+  \ib {f''_i}(\partial_i u) |\partial_i\partial_i u|^2 \Gamma^{\beta_i}(|\partial_i u|) \eta^{2k}\dx + I_{2,i}\Bigg] \, . 
\end{eqnarray}

Let us consider $I_{2,i}$, $i\in \{1, \dots ,n\}$. With $\beta_i > - 1/2$ as above we have
\begin{eqnarray}\label{proof 12}
I_{2,i} &=& \ib |\partial_i\partial_i u| {f''_i}^{\frac{1}{2}}(\partial_i u) \Gamma^{\frac{\beta_i}{2}}(|\partial_i u)
(f''_i)^{-\frac{1}{2}}(\partial_i u) \Gamma^{-\frac{\beta_i}{2}}(|\partial_i u|)\nonumber\\[2ex]
&& \quad \cdot \Gamma^{\frac{1}{2}+\gamma_{i}}(|\partial_i u|)  |f_i'|(\partial_i u)\eta^{2k}\dx\nonumber\\[2ex]
&\leq & c \ib {f''_i} (\partial_i u) |\partial_i \partial_i u|^2 \Gamma^{\beta_i}(|\partial_i u|) \eta^{2k}\dx \nonumber\\[2ex]
&& + c \ib (f''_i)^{-1}(\partial_i u) 
\Gamma^{1+\gamma_{i}-\eps_i}(|\partial_i u|) |f_i'|^2(\partial_i u)  \eta^{2k}\dx \, .
\end{eqnarray}
The first integral on the right-hand side of \reff{proof 12} already occurs in \reff{proof 11} and 
the second one is handled with \reff{pre 2} and Lemma \ref{proof lem 2} (recalling $\eps_i > \theta_i$)
\begin{eqnarray}\label{proof 13}
\lefteqn{\ib (f''_i)^{-1}(\partial_i u) \Gamma^{1+\gamma_{i}-\eps_i}(|\partial_i u|) |f_i'|^2(\partial_i u)  \eta^{2k}\dx}\nonumber\\[2ex]
& \leq & \ib f_i (\partial_i u) \Gamma^{1+\gamma_{i}-(\eps_i-\theta_i)}(|\partial_i u|)\eta^{2k} \dx\nonumber\\[2ex]
 & \leq & \ib \Big[f_i (\partial_i u) \Gamma^{1+\gamma_{i}}(|\partial_i u|)\Big]^{\frac{1}{\rho}}\eta^{\frac{2k}{\rho}}
\eta^{\frac{2k}{\rho*}}\dx + c \nonumber\\[2ex]
 &\leq & \eps \ib f_i (\partial_i u) \Gamma^{1+\gamma_{i}}(|\partial_i u|)\eta^{2k}\dx + c(\eps,r) 
\end{eqnarray}
and once more the integral on the right-hand side is absorbed.\\

To sum up, \reff{proof 11} implies with the help of \reff{proof 12} and \reff{proof 13} for $i=1$, \dots ,$n$
\begin{eqnarray}\label{proof 14}
\lefteqn{\ib f_i (\partial_i u) \Gamma^{1+\gamma_{i}}(|\partial_i u|) \eta^{2k}\dx}\nonumber\\[2ex]
& \leq &
c \Bigg[1+  \ib f''_i(\partial_i u)  |\partial_i\partial_i u|^2 \Gamma^{\beta_i}(|\partial_i u|) \eta^{2k}\dx\Bigg] \, .
\end{eqnarray}

Discussing the right-hand side of \reff{proof 14} we apply Lemma \ref{proof lem 1}, where we let $f(Z) = \sum_{j=1}^n f_j(Z_j)$
and fix $i\in\{1,\dots , n\}$:
\begin{eqnarray}\label{proof 15}
\lefteqn{\ib f''_i(\partial_i u) |\partial_i\partial_i u| ^2 \Gamma^{\beta_i}(|\partial_i u|) \eta^{2k}\dx}\nonumber\\[2ex]
 &\leq &
c \ib D^2 f(\nabla u)\big(\partial_i \nabla u, \partial_i \nabla u\big) \Gamma^{\beta_i}(|\partial_i u|) \eta^{2k}\dx \nonumber\\[2ex]
&\leq & c \ib D^2 f(\nabla u) \big(\nabla \eta,\nabla \eta)  \Gamma^{1+\beta_i}(|\partial_i u|)\eta^{2k-2}\dx \nonumber\\[2ex]
&\leq & c(r) \sum_{j=1}^n \ib f''_j(\partial_j u)  \Gamma^{1+\beta_i}(|\partial_i u|) \eta^{2k-2}\dx\, .
\end{eqnarray}
For $j=i$ on the right-hand side of \reff{proof 15} we now apply the left-hand side of \reff{pre 1}
and again \mbox{Lemma \ref{proof lem 2}} with the result (recall $\delta_i + \eps_i <1$)
\begin{eqnarray}\label{proof 16}
\lefteqn{ \ib f''_i(\partial_i u) \Gamma^{1+\beta_i}(|\partial_i u|) \eta^{2k-2}\dx}\nonumber\\[2ex]
& \leq & \ib f_i(\partial_i u) \Gamma^{\gamma_{i}+\delta_i+\eps_i}(|\partial_i u|)\eta^{2k-2}\dx \nonumber\\[2ex]
&\leq &  \ib \Big[f_i(\partial_i u) \Gamma^{1+\gamma_{i}}(|\partial_i u|))\Big]^{\frac{1}{\rho}} \eta^{\frac{2k}{\rho}}
 \eta^{\frac{2k}{\rho^*}-2}\dx +c\nonumber\\[2ex]
 &\leq & \eps \ib f_i(\partial_i u) \Gamma^{1+\gamma_{i}}(|\partial_i u|)\eta^{2k}\dx + c(\eps,r) \, .
\end{eqnarray}
Note that the integral on the right-hand side of \reff{proof 16} can be absorbed in the left-hand side of
\reff{proof 14}. This proves Proposition \ref{proof prop 2}. \qed\\

\section{Iteration}\label{iter}

We start with an elementary proposition recalling and relating the relevant parameters of the problem.

\begin{proposition}\label{iter prop 1} With $q_i^\pm$, $\underline{q}_i$, $\overline{q}_i$,
$\delta_i$, $\theta_i$, $\beta_i$, $\gamma_i$, $\eps_i$, $i=1$, \dots , $n$
as above we further let $\vart =: 1-\delta_i$ and
\[
\omega_i^\pm := \frac{q_i^{\pm}}{2} + \gamma_i \, ,  \quad i \in\{1,\dots , n\}\, .
\]
We fix $\tau \geq 0$, $i$, $j \in \{1,\dots ,n\}$  and choose $\gamma_i$ such that
($M >0$ denoting  an arbitrary fixed number)
\begin{equation}\label{iter 1}
1 + \gamma_i < \left\{\begin{array}{ccl}
\dis \frac{\underline{q}_i\vart}{2} \frac{2+\tau}{\overline{q}_j -2\vart}  - \eps_i \vart\frac{\tau+\overline{q}_j/\vart}{\overline{q}_j -2\vart}
& \mbox{if}& \overline{q}_j >2\vart \\[3ex]
\dis M& \mbox{if} & \overline{q}_j \leq 2  \vart\end{array}\right\} \, .
\end{equation}
This yields (for any combination of $q_j^{\pm}$ and $q_i^{\pm}$)
\begin{equation}\label{iter 2}
q^{\pm}_j \frac{1+\beta_i}{\omega_i^{\pm} - \beta_i} < 
2 \vart \frac{1+\frac{q_i^{\pm}}{2} + \gamma_i}{\omega_i^{\pm}-\beta_i} +
\tau\vart  \, ,
\end{equation}
\end{proposition}

\emph{Proof.} In the case $\overline{q}_j >2 \vart$ we note that
\[
1+ \gamma_i < \frac{\underline{q}_i\vart}{2}\frac{2+\tau}{\overline{q}_j-2\vart}
 - \eps_i \vart\frac{\tau+\overline{q}_j/\vart}{\overline{q}_j -2\vart} \, ,
\]
which is equivalent to
\[
(1+ \gamma_i) \big[\overline{q}_j - 2 \vart\big] < 
\underline{q}_i\vart + \tau \Bigg[\frac{\underline{q}_i\vart}{2}-\eps_i\vart\Bigg]  - \eps_i \overline{q}_j\, .
\]
Writing this in the form
\[
\overline{q}_j(1+\beta_i) < 2 \vart\Bigg[1+ \gamma_i + \frac{\underline{q}_i}{2}\Bigg]+ 
\tau \vart\Bigg[\frac{\underline{q}_i}{2}-\eps_i\Bigg]
\]
and recalling that we have by definition $\omega^{\pm}_i - \beta_i = (q_i^\pm/2) - \eps_i$ we obtain as an equivalent inequality
\[
\overline{q}_j \frac{1+\beta_i}{\omega_i^{\pm} - \beta_i} < 
2 \vart \frac{1+\frac{\underline{q}_i}{2} + \gamma_i}{\omega_i^{\pm}-\beta_i} +
\tau \vart \frac{\underline{q}_i-2\eps_i}{q_i^{\pm}-2\eps_i} \, . \qed
\]

\vspace*{2ex}
Up to now no relation between $q_i^+$ and $q_i^-$ was needed due to our particular Ansatz depending on $t$ instead of $|t|$.\\

To complete the proof of Theorem \ref{main 2} it remains to handle the mixed terms on the right-hand side of \reff{proof 7}. 
Here, of course, it is no longer possible to argue with the structure conditions for fixed $i$, i.e. to 
argue with $q_i^\pm$ separated from each other  in disjoint regions.\\

Throughout the rest of this section we suppose that the assumptions of Theorem \ref{main 2} are satisfied.\\

Consider a set $U \subset \Omega$ and a $C^1$-function $v$: $\Omega \to \rz$. 
We let for any $i \in \{1, \dots , n\}$
\[
U \cap [\partial_i v \geq 0] =: U_i^+[v] =: U_i^+\, , \qquad 
U \cap [\partial_i v < 0] =: U_i^-[v] =: U_i^-\, , 
\]
in particular $u$ can be written as the disjoint union
\[
U = U_i^+ \cup U_i^- \, .
\]

and for every $1 \leq i \leq n$.\\

Using this notation, recalling Proposition \ref{proof prop 2} and the left-hand side of \reff{pre 1} we have for every $1 \leq i \leq n$
\begin{eqnarray}\label{iter 3}
\lefteqn{\ib f_i (\partial_i u) \Gamma^{1+\gamma_{i}}(|\partial_i u|)  \eta^{2k}\dx}\nonumber\\[2ex] 
&\leq &c \sum_{j\not= i}  \ib f''_{j}(\partial_j u) \Gamma^{1+\beta_i}(|\partial_i u|) \eta^{2k-2}\dx\nonumber\\[2ex]
&\leq & c \sum_{j\not= i}\ib f_j(\partial_j u) \Gamma^{\delta_i -1}(|\partial_j u|) \Gamma^{1+\beta_i}(|\partial_i u|)\eta^{2k-2}\dx 
\end{eqnarray}

Fix $i \in \{1, \dots , n\}$. For any $j\in \{1, \dots , n\}$ we let
\[
\kappa_i^{\pm} = \frac{1+\omega_i^\pm}{1+\beta_i}\, ,\qquad
\hat{\kappa}_i^{\pm} = \frac{1+\omega_i^{\pm}}{\omega_i^{\pm} - \beta_i}\, .
\]
This gives for fixed $1 \leq i \leq n$ and for $\eps > 0$ sufficiently small (note that the ball $B$ is divided into two parts w.r.t.~the
function $\partial_i u$)
\begin{eqnarray}\label{iter 4}
\lefteqn{\sum_{j\not= i}\ib   f_j(\partial_j u) \Gamma^{\delta_i-1}(|\partial_j u|)\Gamma^{1+\beta_i}(|\partial_i u|) \eta^{2k-2}\dx}\nonumber\\[2ex]
&\leq c & \sum_{j\not=i}\sum_{\pm} \int_{B_{i,\pm}} \big(1+f_j(\partial_j u)\big)
 \Gamma^{\delta_i -1}(|\partial_j u|)  \Gamma^{1+\beta_i}(|\partial_i u|) 
\eta^{2k-2}\dx\nonumber\\[2ex]
&\leq&  \sum_{j\not= i}\sum_{\pm} \Bigg[\eps \int_{B_{i,\pm}} \Gamma(|\partial_i u|)^{1+\omega_i^\pm} \eta^{2k}\dx \nonumber\\
&& + c(\eps) \int_{B_{i,\pm}} \big(1+f_j(\partial_j u)\big)^{\frac{1+\omega_i^\pm}{\omega_i^{\pm}-\beta_i}} 
\Gamma^{(\delta_i -1)\frac{1+\omega_i^\pm}{\omega_i^\pm - \beta_i}}(|\partial_j u|)\dx \Bigg] \, .
\end{eqnarray}
By \reff{pre 3} we have on $B_{i,\pm}$ for $|\partial_i u|$ sufficiently large 
$\Gamma (|\partial_i u|)^{q_i^{\pm}/2} \leq c f_i(\partial_i u)$, hence
by the definition of $\omega_i^{\pm}$
\begin{eqnarray*}
\lefteqn{\eps \sum_{j\not= i}\sum_{\pm} \int_{B_{i,\pm}} \Gamma(|\partial_i u|)^{1+\omega_i^\pm} \eta^{2k}\dx}\\[2ex] 
&& \leq (n-1) \eps \ib f_i(\partial_i u) \Gamma^{1+\gamma_i}(|\partial_i u|) \eta^{2k}\dx + c
\end{eqnarray*}
and, as usual, the integral on the right-hand side can be absorbed in \reff{iter 3}.\\

We will finally show with the help of an iteration procedure that for every $1 \leq i \leq n$
\begin{equation}\label{iter 5}
\sum_{j\not= i} \sum_{\pm}\int_{B_{i,\pm}} \big(1+f_j(\partial_j u)\big)^{\frac{1+\omega_i^\pm}{\omega_i^{\pm}-\beta_i}} 
\Gamma^{(\delta_i -1) \frac{1+\omega_i^\pm}{\omega_i^\pm - \beta_i}}(|\partial_j u|) \dx \leq c\, ,
\end{equation}
which completes the proof of Theorem \ref{main}.\\

If fact, let us suppose that  \reff{iter 1} is true with the choice $\vart = 1 - \delta_i$.
Then we may apply Proposition \ref{iter prop 1} and \reff{iter 2} implies 
in the case $\overline{q}> 2 \vart$
\[
\Gamma^{(\delta_i -1) \frac{1+\omega_i^\pm}{\omega_i^\pm - \beta_i} + (\delta_i-1) \frac{\tau }{2} }(|\partial_j u|)
\leq c \big(1+f_j(\partial_j u)\big)^{-\frac{1+\beta_i}{\omega_i^\pm - \beta_i}} \, .
\]
Thus we obtain
\begin{equation}\label{iter 6}
\big(1+f_j(\partial_j u)\big)^{\frac{1+\omega_i^\pm}{\omega_i^{\pm}-\beta_i}} 
\Gamma^{(\delta_i -1) \frac{1+\omega_i^\pm}{\omega_i^\pm - \beta_i}}(|\partial_j u|)
\leq c  \big(1+f_j(\partial_j u)\big) \Gamma^{(1-\delta_i) \frac{\tau}{2}}(|\partial_j u|) \, .
\end{equation}
In the case $\overline{q}_j \leq 2 \vart$ we have
\[-
(1+f_j(\partial_j u) \leq c \Gamma^{1-\delta_i}(|\partial_j u|) \, ,
\]
hence
\[
\big(1+f_j(\partial_j u)\big)^{\frac{1+\omega_i^\pm}{\omega_i^{\pm}-\beta_i}} 
\Gamma^{(\delta_i - 1) \frac{1+\omega_i^\pm}{\omega_i^\pm - \beta_i}}(|\partial_j u|)
\leq c 
\]
and \reff{iter 6} holds as well.\\ 

We note that \reff{iter 6} is formulated uniformly w.r.t.~the index $j$ and
the symbol $\pm$ is just related to $\partial_i u$.\\

Inequality \reff{iter 6} is the main tool for the following iteration leading to the claim \reff{iter 5}.\\

\underline{$i=1$.}\\

Choosing $\gamma_i > - 1/2 + \theta_i$ sufficiently close to $- 1/2 + \theta_i$, \reff{iter 1} is valid 
with the choice $\tau = 0$ if we have

\begin{equation}\label{iter 7}
\overline{q}_j < \frac{2 \underline{q}_i (1-\delta_i)}{1+2 \theta_i} + 2 (1-\delta_i)  
\qquad \mbox{for all}\qquad 2 \leq  j \leq n\, ,
\end{equation}
and \reff{iter 7} is just assumption \reff{pre 4} for $i=1$.\\

From \reff{iter 1} we deduce \reff{iter 6} for $i=1$ and  \reff{iter 5} follows from \reff{iter 6} 
for $i=1$ and for all $2 \leq j \leq n$ with the choice $\tau =0$
\begin{eqnarray}\label{iter 8}
\lefteqn{\sum_{j\not= i}\sum_{\pm}\int_{B_i,\pm}\big(1+f_j(\partial_j u)\big)^{\frac{1+\omega_1^\pm}{\omega_1^{\pm}-\beta_1}} 
\Gamma^{(\delta_i -1) \frac{1+\omega_1^\pm}{\omega_1^\pm - \beta_1}}(|\partial_j u|)}\nonumber\\[2ex]
&&\leq  c  \ib \big(1+f_j(\partial_j u)\big) \dx \leq c\, .
\end{eqnarray}

Returning to \reff{iter 3} and \reff{iter 4} we insert \reff{iter 8} and 
on account of $1+\gamma_i > 1/2$ we have
\begin{equation}\label{iter 9}
\ib f_1(\partial_1 u) \Gamma^{\frac{1}{2}}(|\partial_1 u|) \eta^{2k}\dx \leq c \, .
\end{equation}

\begin{remark}\label{proof rem 2}
In \cite{BFZ:2007_1} we have $\delta_i=\theta_i =0$, $i=1$, $2$, and w.l.o.g.~the case $p= q_2 \leq q_1 = q$ is considered.
Moreover, $\overline{q}_j = \underline{q}_j = q_j$, $j=1$, $2$.
In this case we trivially have \reff{iter 7}.\\
\end{remark}

\underline{$1 < i \leq n$.}\\

Suppose that we have in addition to \reff{iter 7} (again compare \reff{pre 4})
\begin{equation}\label{iter 10}
\overline{q}_j < \overline{q}_j < \frac{2 \underline{q}_i (1-\delta_i)}{1+2 \theta_i} + 2 (1-\delta_i)  
\qquad\mbox{for}\qquad i+1 \leq j \leq n \, .
\end{equation}
With the same argument leading to \reff{iter 8} we have for all $i+1 \leq j \leq n$
\begin{equation}\label{iter 11}
\sum_{j> i} \sum_{\pm}\int_{B_i,\pm}\big(1+f_j(\partial_j u)\big)^{\frac{1+\omega_i^\pm}{\omega_i^{\pm}-\beta_i}} 
\Gamma^{(\delta_i -1) \frac{1+\omega_i^\pm}{\omega_i^\pm - \beta_i}}(|\partial_j u|)
\leq c  \, .
\end{equation}

Moreover, we suppose that by iteration we have \reff{iter 9} for $1 \leq j < i$, i.e.
\begin{equation}\label{iter 12}
\ib f_j(\partial_j u) \Gamma^{\frac{1}{2}}(|\partial_j u|) \eta^{2k}\dx \leq c \,, \quad 1 \leq j < i \, .
\end{equation}

Then we return to \reff{iter 1} with the choice $\tau = (1-\delta_i)^{-1}$. 
For $\gamma_i > -1/2+\theta_i$ and $\gamma_i$
sufficiently close to $-1/2+\theta_i$ we are lead to the condition
\begin{equation}\label{iter 13}
\overline{q}_j < \frac{2}{1+2\theta_i}\Bigg[\frac{\underline{q}_i}{2}(1-\delta_i)\Big[2+\frac{1}{1- \delta_i}\Big] 
- \theta_i (1+\overline{q}_j)\Bigg]+2 (1-\delta_i)\, ,
\end{equation}
$1 \leq j < i$, and \reff{iter 13} is just the assumption \reff{pre 5}.\\

With \reff{iter 1} we again have \reff{iter 6}, now with $\tau = (1-\delta_i)^{-1}$, hence
\begin{eqnarray}\label{iter 14}
\lefteqn{\sum_{j < i}\sum_{\pm} \int_{B_i,\pm} \big(1+f_j(\partial_j u)\big)^{\frac{1+\omega_i^\pm}{\omega_i^{\pm}-\beta_i}} 
\Gamma^{(\delta_i -1)\frac{1+\omega_i^\pm}{\omega_i^\pm - \beta_i}}(|\partial_j u|)\dx}\nonumber\\[2ex]
&& \leq   c \sum_{j< i} \ib \big(1+f_j(\partial_j u)\big) \Gamma^{\frac{1}{2}}(|\partial_j u|)\dx \leq c  \, ,
\end{eqnarray}

With \reff{iter 11} and \reff{iter 14} one has
\begin{equation}\label{iter 15}
\sum_{j\not= i}\sum_{\pm}\int_{B_i,\pm}\big(1+f_j(\partial_j u)\big)^{\frac{1+\omega_1^\pm}{\omega_1^{\pm}-\beta_1}} 
\Gamma^{(\delta_i -1)\frac{1+\omega_1^\pm}{\omega_1^\pm - \beta_1}}(|\partial_j u|) \leq c\, ,
\end{equation}
which exactly as in the case $i=1$ shows
\begin{equation}\label{iter 16}
\ib f_i(\partial_i u) \Gamma^{\frac{1}{2}}(|\partial_i u|) \eta^{2k}\dx \leq c \, ,
\end{equation}
hence with \reff{iter 16} we proceed one step in the iteration of \reff{iter 13}.
This completes the proof of \mbox{Theorem \ref{main}}. \qed\\

\begin{remark}\label{proof rem 3}
With the notation of Remark \ref{proof rem 2} condition \reff{iter 1} reduces in the situation discussed in \cite{BFZ:2007_1} to
\[
1 + \gamma_1 < \frac{q}{2} \frac{3}{p-2}
\]
on account of
\[
\frac{q}{2}\frac{3}{p-2} \geq \frac{q}{2}\frac{3}{q-2} > \frac{3}{2}
\]
we may choose $\gamma_1 =1/2$ without imposing any condition relating $p$ and $q$, thus the results presented in
\cite{BFZ:2007_1} immediately follow as a corollary. 
\end{remark}





\bibliography{asymmetrical_splitting}
\bibliographystyle{unsrt}

\end{document}